\documentclass[letterpaper, 10 pt, conference]{ieeeconf}  

\IEEEoverridecommandlockouts                              

\usepackage{graphics} 
\usepackage{epsfig} 
\usepackage{mathptmx} 
\usepackage{times} 
\usepackage{amsmath} 
\usepackage{amssymb}  
 
\usepackage{amsthm}
\usepackage{lipsum}
\usepackage{stix}
\usepackage{tikz}
\usetikzlibrary{patterns}
\usepackage{multirow}

\newtheorem{remark}{Remark}
\newtheorem{problem}{Problem}
\newtheorem{definition}{Definition}

\newtheorem{example}{Example}

\title{\LARGE \bf
Orientation Control of the Bouncing Ball*
}

\author{William Clark$^{1}$ and Dora Kassabova$^{1}$
\thanks{*This work was funded by NSF grant DMS-1645643.}
\thanks{$^{1}$William Clark and Dora Kassabova are with the Department of Mathematics, Cornell University, Ithaca, NY 14850, USA
        {\tt\small wac76@cornell.edu}, {\tt\small dmk285@cornell.edu}}%
}

\begin{document}

\maketitle
\thispagestyle{empty}
\pagestyle{empty}

\begin{abstract}

Control of a hybrid dynamical system can manifest in one of two main ways: either through the continuous or the discrete dynamics. An example of controls influencing the continuous dynamics is legged locomotion, where the joints are actuated but the location and nature of the impacts are uncontrolled. In contrast, an example of discrete control would be in tennis; the player can only influence the trajectory of the ball through striking it.

This work examines the latter case with two key emphases. The first is that controls manifest through changing the location of the guard (as opposed to changing only the reset). The second is that the location of the guard is described by ``external variables'' while the goal is to control ``internal variables.'' As a simple test of this theory, orientation control of a bouncing ball is explored; the ball is only controlled during impacts which are exclusively position-dependent.

\end{abstract}

\section{INTRODUCTION}\label{sec:intro}
Systems which undergo both continuous- and discrete-time transitions are referred to as hybrid. The dynamics of such a system are described by
\begin{equation*}
    \begin{cases}
        \dot{x} = f(x), & x\in M\setminus S, \\
        x^+ = \Delta(x^-), & x^-\in S,
    \end{cases}
\end{equation*}
where
\begin{enumerate}
    \item $M$ is a (finite-dimensional) manifold,
    \item $S\subset M$ is an embedded codimension 1 submanifold,
    \item $f\in\mathfrak{X}(M)$ is a vector-field, and
    \item $\Delta:S\to M$ is a map.
\end{enumerate}
All the data $(M,S,f,\Delta)$ is assumed to be smooth. The ambient manifold, $M$, will be referred to as the state-space, $S$ is the guard, and $\Delta$ is the reset.

Controls can be added to this system in two main ways. If controls are added to the continuous component, the resulting system would have the form
\begin{equation*}
    \begin{cases}
        \dot{x} = f(x,u), & x\in M\setminus S, \\
        x^+ = \Delta(x^-), & x^-\in S.
    \end{cases}
\end{equation*}
A major application of this type of control problem is in legged locomotion \cite{HBG_walking_2016,RL_top_2021,HG_2016,HG_2019,GK_2017}. On the other hand, controls can be inserted into the reset conditions,
\begin{equation}\label{eq:reset_control_hybrid}
    \begin{cases}
        \dot{x} = f(x), & x\in M\setminus S_u, \\
        x^+ = \Delta(x^-,u), & x^-\in S_u.
    \end{cases}
\end{equation}
A class of controlled hybrid systems falling under this type are juggling systems, \cite{BRK1999,RK19925,kant2022,kant2021,lynch2001,gerard2004,ronsse2006} as well as systems with impulsive controls \cite{kmk2021}.

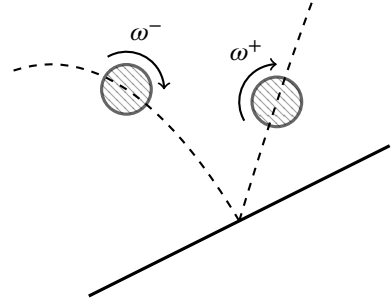
\begin{figure}
    \centering
    \begin{tikzpicture}
        \draw [very thick, black] (0,0) -- (4,2);
        \draw [thick, dashed, domain=-1:2] plot (\x, {3-1/3*\x-1/3*\x*\x});
        \draw [thick, dashed, domain=2:3] plot (\x, {-7 + 14/3*\x-1/3*\x*\x});
        \fill [very thick, pattern=north west lines,opacity=.6,draw] (0.5,2.75) circle [radius=0.33];
        \fill [very thick, pattern=north west lines,opacity=.6,draw] (2.5,2.5833) circle [radius=0.33];
        \draw[thick, ->] (0.25,3.183) arc (120:0:0.5);
        \draw[thick, ->] (2.0670,2.3333) arc (210:90:0.5);
        \node at (0.7736,3.5018) {$\omega^-$};
        \node at (2.1,3.2761) {$\omega^+$};
    \end{tikzpicture}
    \caption{As a ball strikes and ricochets off a surface, its angular and linear velocity will generally interact. As such, it seems natural that the pre- and post-impact angluar velocities will differ, i.e. $\omega^-\ne\omega^+$.}
    \label{fig:intro_picture}
\end{figure}

An important subclass of hybrid systems are those of impact type, \cite{brogliato2016,fetecau2003}, and will be the focus of this work. 
The typical setup of an impact system is a hybrid system with the additional structure of a natural Lagrangian system, and can be generated by the following three pieces of data $(Q,L,h)$ where
\begin{enumerate}
    \item $(Q,g)$ is a (finite-dimensional) Riemannian manifold.
    \item $L:TQ\to\mathbb{R}$ is a natural Lagrangian,
    \begin{equation*}
        L(v) = \frac{1}{2}g(v,v) - \tau_Q^*V(v),
    \end{equation*}
    where $V:Q\to\mathbb{R}$ is the potential and $\tau_Q:TQ\to Q$ is the canonical tangent bundle projection.
    \item $h:Q\to \mathbb{R}$ is a smooth function with zero as a regular value.
\end{enumerate}
An impact system generates a hybrid system via $M=TQ$, $f=X_L$ is the Lagrangian vector-field, $S$ is given by outward pointing vectors, i.e.
\begin{equation}\label{eq:guard}
    S = \left\{ v\in TQ : \tau_Q^*h(v) = 0, \quad dh(v) > 0 \right\},
\end{equation}
and the reset map (assuming elastic impacts) is given by the Weierstrass-Erdmann conditions
\begin{equation}\label{eq:still_WE_condition}
    \begin{split}
        \frac{\partial L}{\partial v}^+ - \frac{\partial L}{\partial v}^- &= \varepsilon\cdot dh, \\
        -E_L^+ + E_L^- &= 0,
    \end{split}
\end{equation}
where $\varepsilon$ is chosen to ensure conservation of energy, $E_L$. When the Lagrangian is natural, these conditions specialize to
\begin{equation}\label{eq:uncontrolled_elastic_impact}
    \Delta(v) = v - 2\frac{dh(v)}{dh(\nabla h)} \nabla h,
\end{equation}
where $g(\nabla h, \cdot) = dh$. 

Controls will be implemented in an impact system by changing the event function to depend on the following controls $h:\mathcal{U}\times Q\to\mathbb{R}$, where the set $\mathcal{U}$ contains all admissible controls. In this more general setting, \eqref{eq:uncontrolled_elastic_impact} is no longer valid and the correct reset conditions are given by
\begin{equation}\label{eq:moving_WE_condition}
    \begin{split}
        \frac{\partial L}{\partial v}^+ - \frac{\partial L}{\partial v}^- &= \varepsilon \cdot d_xh, \\
        -E_L^+ + E_L^- &= \varepsilon\cdot \frac{\partial h}{\partial t}.
    \end{split}
\end{equation}
Although \eqref{eq:still_WE_condition} and \eqref{eq:moving_WE_condition} provide the ``correct'' answer via variations, they lack the key qualitative property of symmetry breaking. Suppose that both $L$ and $h$ are invariant under the action of a Lie group $G$; such a system is called a hybrid system with symmetry \cite{amesThesis}. Then the reset map induced by \eqref{eq:still_WE_condition} is symmetry-preserving in the sense that it preserves the momentum map, $\Delta^*J = J|_S$, for which we define
\begin{gather*}
    J:TQ\to\mathfrak{g}^* \\
    J(v)(\xi) = g\left( v, \xi_Q(\tau_Q(v)) \right),
\end{gather*}
where $\mathfrak{g}$ is the Lie algebra of $G$ and $\xi_Q$ is the vector-field on $Q$ generated by the element $\xi\in\mathfrak{g}$. 

A spherical ball has the symmetry group $\mathrm{SO}_3$ and, by the discussion above, its angular momentum is conserved across impacts. A reason why this property is too simplistic is that it would be impossible to influence its angular velocity; in the situation shown in Figure \ref{fig:intro_picture} we would always have the equality $\omega^-=\omega^+$. If this were the case, tennis players would be unable to control the spin of the tennis ball. 

The goal of this work is two-fold. First, we present a modified version of \eqref{eq:still_WE_condition} and \eqref{eq:moving_WE_condition} such that the momentum is no longer conserved after the reset. This is accomplished by enforcing a nonholonomic constraint at the moment of impact; the ball must roll without slipping while in contact with the surface. Second, we use this impact condition to study the controllability of the orientation of a planar bouncing disk, which is an impact hybrid system of the form \eqref{eq:reset_control_hybrid}.

The classical case of impact systems with symmetries is reviewed in \S\ref{sec:symmetries}. By imposing nonholonomic constraints at the moment of impact, the symmetry breaking impact law is constructed in \S\ref{sec:break_impacts}. Symmetry breaking impacts are reformulated into the intermittent control problem, \eqref{eq:reset_control_hybrid}, in \S\ref{sec:control_problem}. To test this theory, the example of orientation control of a planar bouncing ball is explored in \S\ref{sec:bouncing_ball}. Conclusions and future work is in \S\ref{sec:conclusions}.

\section{IMPACT SYSTEMS WITH SYMMETRIES}\label{sec:symmetries}
There have been works done in dealing with symmetries inside impact systems \cite{amesThesis} and recently in \cite{colombo2020} (and the references therein). We first introduce the classical continuous case and later extend to include impacts.
\subsection{Continuous Systems with Symmetries}
All Lagrangian functions here will be assumed to be of natural type.
\begin{definition}
    For a smooth (finite-dimensional) manifold $Q$, a natural Lagrangian is a function $L:TQ\to\mathbb{R}$ of the form
    \begin{equation*}
        L(v) = \frac{1}{2}g(v,v) - \tau^*V(v),
    \end{equation*}
    where $g$ is a Riemannian metric on $Q$, $V:Q\to\mathbb{R}$ is a smooth function and $\tau:TQ\to Q$ is the canonical tangent bundle projection.
\end{definition}
A Lagrangian system with symmetry will be defined as follows:
\begin{definition}
    A Lagrangian system with symmetry is given by the 4-tuple $(Q,L,G,\Phi)$ where
    \begin{enumerate}
        \item $L$ is a natural Lagrangian,
        \item $G$ is a Lie group with a free and proper action $\Phi:G\times Q\to Q$,
        \item the Riemannian metric and potential function are invariant under the group action.
    \end{enumerate}
\end{definition}

As a consequence of Noether's theorem, the induced momentum map is preserved, $J:TQ\to\mathfrak{g}^*$,
\begin{equation*}
    J(v)(\xi) = g\left( v, \xi_Q(\tau(v))\right),
\end{equation*}
where $\mathfrak{g}^*$ is the dual to the Lie algebra of $G$. By using the locked intertia tensor, $\mathbb{I}_x:\mathfrak{g}\to\mathfrak{g}^*$, the momentum map can be turned into the mechanical connection. The locked inertia tensor is given by
\begin{equation*}
    \mathbb{I}_x(\xi)(\zeta) = g\left( \xi_Q(x),\zeta_Q(x)\right), \quad \xi,\zeta\in\mathfrak{g},
\end{equation*}
where $\xi_Q\in\mathfrak{X}(Q)$ is the vector-field induced by differentiating the group action in the direction $\xi\in\mathfrak{g}$.
With this, we define the mechanical connection $\mathcal{A}:TQ\to\mathfrak{g}$ as
\begin{equation*}
    \mathcal{A}(v) = \mathbb{I}^{-1}_{\tau(v)}J(v).
\end{equation*}
This map is, again, preserved under the dynamics.
\begin{example}\label{ex:disk_continuous}
    Consider a disk moving in the plane. Its configuration space is given by $\mathrm{SE}_2=\mathbb{R}^2\times S^1$ and its Lagrangian by
    \begin{equation*}
        L = \frac{1}{2}m\left(\dot{x}^2 + \dot{y}^2\right) + \frac{1}{2}I\dot\theta^2 - V(x,y).
    \end{equation*}
    The symmetry group acting on this system is $G=S^1$ with the action $\Phi(s,(x,y,\theta)) = (x,y,\theta+s)$. Conservation of angular momentum follows from the momentum map,
    \begin{equation*}
        J(x,y,\theta;\dot{x},\dot{y},\dot\theta) = I\dot\theta,
    \end{equation*}
    and, likewise, conservation of angular velocity follows from the connection,
    \begin{equation*}
        \mathcal{A}(x,y,\theta;\dot{x},\dot{y},\dot\theta) = \dot\theta.
    \end{equation*}
\end{example}
\subsection{Impacts with Symmetries}
Here, we expand on symmetries in Lagrangian systems when impacts are present. For unconstrained systems the impact conditions will be chosen such that they are variational on the velocities. If the system is subjected to nonholonomic constraints, the impact conditions will be chosen such that the change in velocity obeys Lagrange-d'Alembert's principle. Under this assumption, impacts have the form
\begin{equation}\label{eq:variational_impact}
    \begin{split}
        \left( \frac{\partial L}{\partial v}^+ - \frac{\partial L}{\partial v}^- \right)\cdot \delta x &= 0 \\
        -\left( E_L^+ - E_L^-\right)\cdot \delta t &= 0.
    \end{split}
\end{equation}
If there are no constraints present and impact occurs when $h(q,t)=0$, the variations satisfy
\begin{equation*}
    \frac{\partial h}{\partial x} \delta x + \frac{\partial h}{\partial t}\delta t = 0,
\end{equation*}
and \eqref{eq:variational_impact} becomes \eqref{eq:moving_WE_condition}. We can now state the definition of a hybrid Lagrangian system with symmetry.
\begin{definition}
    A natural impact Lagrangian system with symmetry is given by the 5-tuple $(Q,L,G,\Phi,h)$ such that
    \begin{enumerate}
        \item $(Q,L,G,\Phi)$ is a natural Lagrangian system with symmetry, and
        \item $h$ is a $G$-invariant function with zero as a regular value.
    \end{enumerate}
\end{definition}
As expected, the induced impact map from a natural impact Lagrangian system with symmetry preserves both the (hybrid) momentum map and mechanical connection, i.e. if $\Delta$ is the impact map generated by \eqref{eq:variational_impact} then
\begin{equation}\label{eq:reset_symmetries}
    \Delta^*J = J|_S, \quad \Delta^*\mathcal{A} = \mathcal{A}|_S,
\end{equation}
where $S$ is the guard given by \eqref{eq:guard}.

\section{SYMMETRY BREAKING IMPACTS}\label{sec:break_impacts}
For a disk moving in the plane, as in Example \ref{ex:disk_continuous}, the angular velocity is always preserved under resets generated by \eqref{eq:variational_impact}. As this would make the orientation of the disk uncontrollable, we seek an alternate formulation of the reset laws. 
If the disk has radius $R>0$ and was constrained to be on the impact surface, its angular and linear velocities are coupled by the rule that it must roll without slipping, i.e. $\dot{x} = R\dot\theta.$

Suppose that we impose $k$ velocity constraints $\eta^1(v)=\eta^2(v)=\ldots=\eta^k(v)=0$ where $\eta^j$ are constraining 1-forms on the ambient space $Q$. Applying these constraints to \eqref{eq:variational_impact}, we end up with the new impact law
\begin{equation}\label{eq:nonh_impact}
    \begin{split}
        \frac{\partial L}{\partial v}^+ - \frac{\partial L}{\partial v}^- &= \varepsilon \cdot d_xh + \lambda_i\cdot \eta^i, \\
        -E_L^+ + E_L^- &= \varepsilon\cdot \frac{\partial h}{\partial t},
    \end{split}
\end{equation}
where the $k+1$ multipliers are chosen to ensure that the constraints are satisfied post-impact, $\eta^i(v) = 0$ for all $i=1,\ldots,k$, see \cite{2101.11128} for more details.

The case of the bouncing ball differs from \eqref{eq:nonh_impact} as the constraints are only present on the guard rather than the ambient space. In particular, the constraining 1-forms $\eta^i$ are on the zero level-set of $h$. For concreteness of notation, let $N = h^{-1}(0)$ be the level-set, so each $\eta^i\in\Omega^1(N)$. By using the metric, these forms can be lifted to the ambient space, $\eta^i_g\in\Omega^1(Q)|_N$. For $x\in N$ and $v\in T_xQ$, let
\begin{equation*}
    \eta^i_g(v) = \eta^i(\pi_g(v)),
\end{equation*}
where $\pi_g:TQ|_N\to TN$ is the $g$-orthogonal projection. We define the modified impact laws as follows.
\begin{definition}
    Let $N\subset Q$ be a co-dimension 1 embedded submanifold with distribution $\mathcal{D}\subset N$. Let $\eta^1,\ldots,\eta^k\in\Omega^1(N)$ be constraints such that $v\in\mathcal{D}$ if and only if $\eta^i(v)=0$ for all $i=1,\ldots,k$. The impact equations are given by
    \begin{equation}\label{eq:new_impact}
    \begin{split}
        \frac{\partial L}{\partial v}^+ - \frac{\partial L}{\partial v}^- &= \varepsilon \cdot d_xh + \lambda_i\cdot \eta^i_g, \\
        -E_L^+ + E_L^- &= \varepsilon\cdot \frac{\partial h}{\partial t}.
    \end{split}
\end{equation}
\end{definition}
\begin{remark}
    The original impact law, \eqref{eq:moving_WE_condition}, assumes that the wall is perfectly slick. By contrast, the new impact law, \eqref{eq:new_impact}, assumes that the wall has an infinite coefficient of friction. Additionally, if the constraint distribution $\mathcal{D}\subset TN$ is not invariant under the group action, the momentum map/connection are no longer preserved, i.e. \eqref{eq:reset_symmetries} is no longer true.
\end{remark}
\begin{remark}
    While there has been work done in studying tangential restitution in impacts, e.g. \cite{domenech2014}, we are unaware of previous attempts to model this phenomenon in a systematic and geometric fashion.
\end{remark}

\section{INTERMITTENT CONTROL PROBLEM} \label{sec:control_problem}
With the symmetry breaking impact law, \eqref{eq:new_impact}, we can almost state the intermittent contact control problem. The control situation is slightly more general than the case presented earlier as the set $N$ is allowed to change as it is influenced by the controls. Consequently, the constraining 1-forms/constraint distribution no longer has a common domain. Below, we list the ingredients for the control problem and later demonstrate a practical way of alleviating the domain issue.

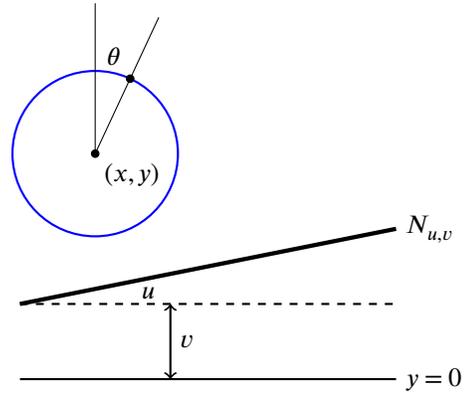
\begin{figure}
    \centering
    \begin{tikzpicture}
        \draw[ultra thick] (0,0) -- (5,1);
        \node[right] at (5,1) {$N_{u,v}$};
        \draw[thick, dashed] (0,0) -- (5,0);
        \node[right] at (1.5,0.15) {$u$};
        \draw[blue, thick] (1,2) circle [radius=1.1];
        \draw[fill=black] (1,2) circle [radius=0.05];
        \node[below right] at (1,2) {$(x,y)$};
        \draw[fill=black] (1.4649,2.9969) circle [radius=0.05];
        \draw (1,2) -- (1,4);
        \draw (1,2) -- (1.8452,3.8126);
        \node[above left] at (1.4649,3.05) {$\theta$};
        \draw[thick] (0,-1) -- (5,-1);
        \node [right] at (5,-1) {$y=0$};
        \draw[thick, <->] (2,-1) -- (2,0);
        \node[right] at (2,-0.5) {$v$};
    \end{tikzpicture}
    \caption{The tabletop is controlled by two parameters; $u$ is its angle of tilt and $v$ is its height.}
    \label{fig:tilted_table}
\end{figure}

The intermittent control problem consists of 
\begin{enumerate}
    \item A Lagrangian system with symmetry, $(Q,L,G,\Phi)$.
    \item A function $h:\mathcal{U}\times Q\to \mathbb{R}$ such that for every admissible control $u\in\mathcal{U}$, zero is a regular-value of $h(u,\cdot):Q\to\mathbb{R}$.
    \item A distribution on each level-set, $\mathcal{D}_u\subset TN_u$ where $N_u = h(u,\cdot)^{-1}(0)$. Moreover, the rank of the distribution is independent of $u$.
\end{enumerate}
Suppose that each distribution is generated by 1-forms in the following sense: $v\in\mathcal{D}_u$ if and only if $\eta_u^i(v)=0$ for all $i$. Call $\iota_u:N_u\hookrightarrow Q$ the inclusion map. Replace each $\eta_u^i$ with a $\tilde{\eta}_u^i\in\Omega(Q)$ such that
\begin{equation*}
    \iota^*_u\tilde{\eta}_u^i = \eta_u^i, \quad
    \tilde{\eta}_u^i|_{N_u} = \eta^i_{u,g}.
\end{equation*}
For regularity purposes, it will be assumed that the forms $\tilde{\eta}_u^i$ depend smoothly on $u$. The controlled version of \eqref{eq:new_impact} becomes
\begin{equation}\label{eq:controlled_impact}
    \begin{split}
        \frac{\partial L}{\partial v}^+ - \frac{\partial L}{\partial v}^- &= \varepsilon\cdot d_xh + \lambda_i\cdot \tilde{\eta}^i_u, \\
        -E_L^+ + E_L^- &= \varepsilon\cdot \frac{\partial h}{\partial u}\cdot \dot{u}.
    \end{split}
\end{equation}
Notice that the location of the impact given by \eqref{eq:controlled_impact} depends on the control and the impact law depends on \textit{both} the control and its time derivative. Energy can be injected or removed from the system by having the control changing at the point of impact.

\begin{figure}
    \centering
    \includegraphics[width=0.45\textwidth]{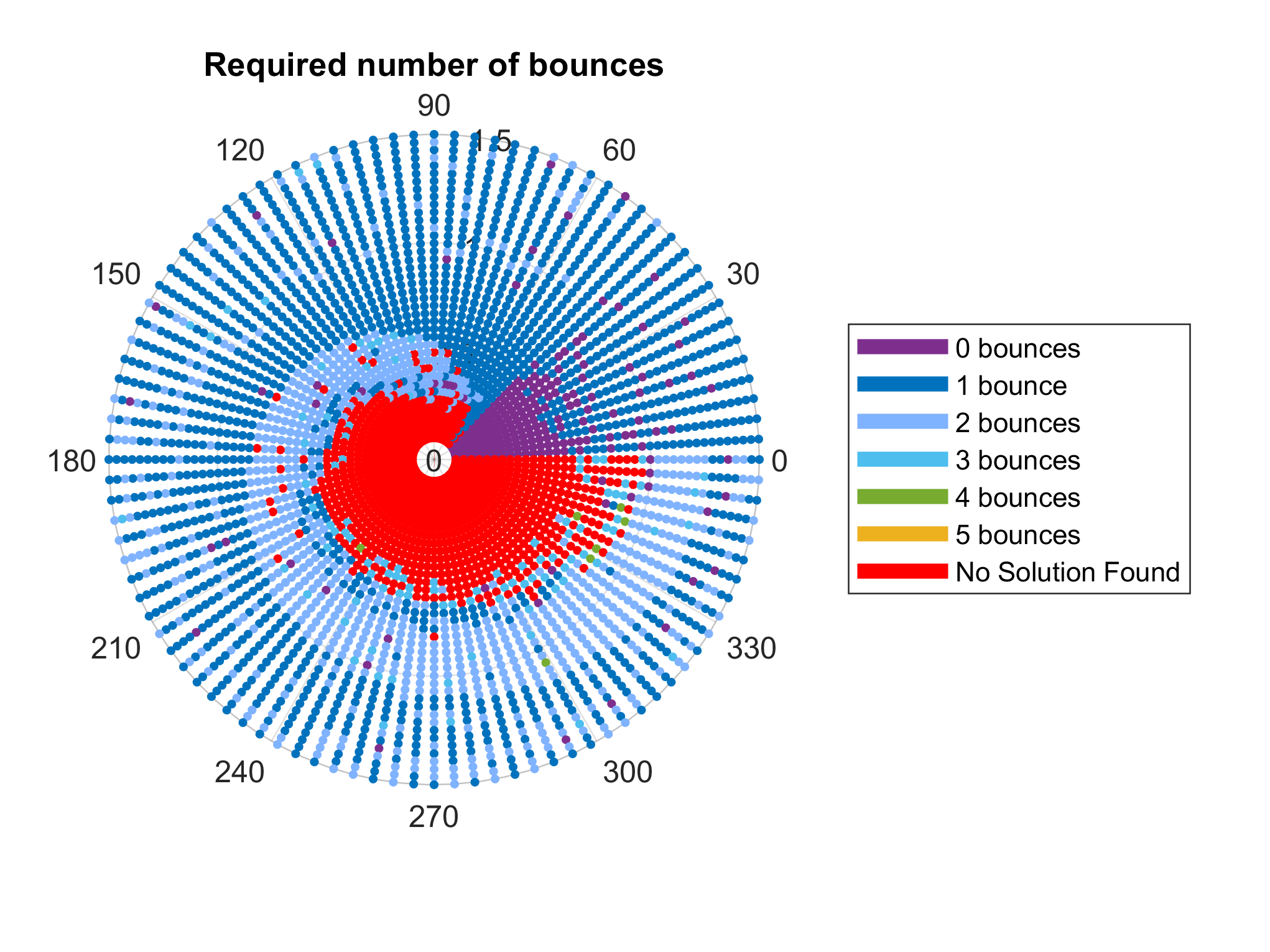}
    \caption{A plot illustrating the orientation controllability of the planar bouncing ball. The radial coordinate indicated total time and the polar coordinate is the desired change in angle.}
    \label{fig:ball_controllability}
\end{figure}

The intermittent control problem can be stated as follows.
\begin{problem}\label{prob:control}
    Let $(Q,L,G,\Phi,\mathcal{U},h,\mathcal{D})$ be an intermittent control problem. For any pair of points $x,y\in Q$ and time $T>0$, does there exist a control law $u:[0,T]\to\mathcal{U}$ that drives the system from $x$ to $y$ while obeying the standard Euler-Lagrange equations away from impacts and subject to the impact law \eqref{eq:controlled_impact}?
\end{problem}

\begin{figure}
    \centering
    \includegraphics[width=0.45\textwidth]{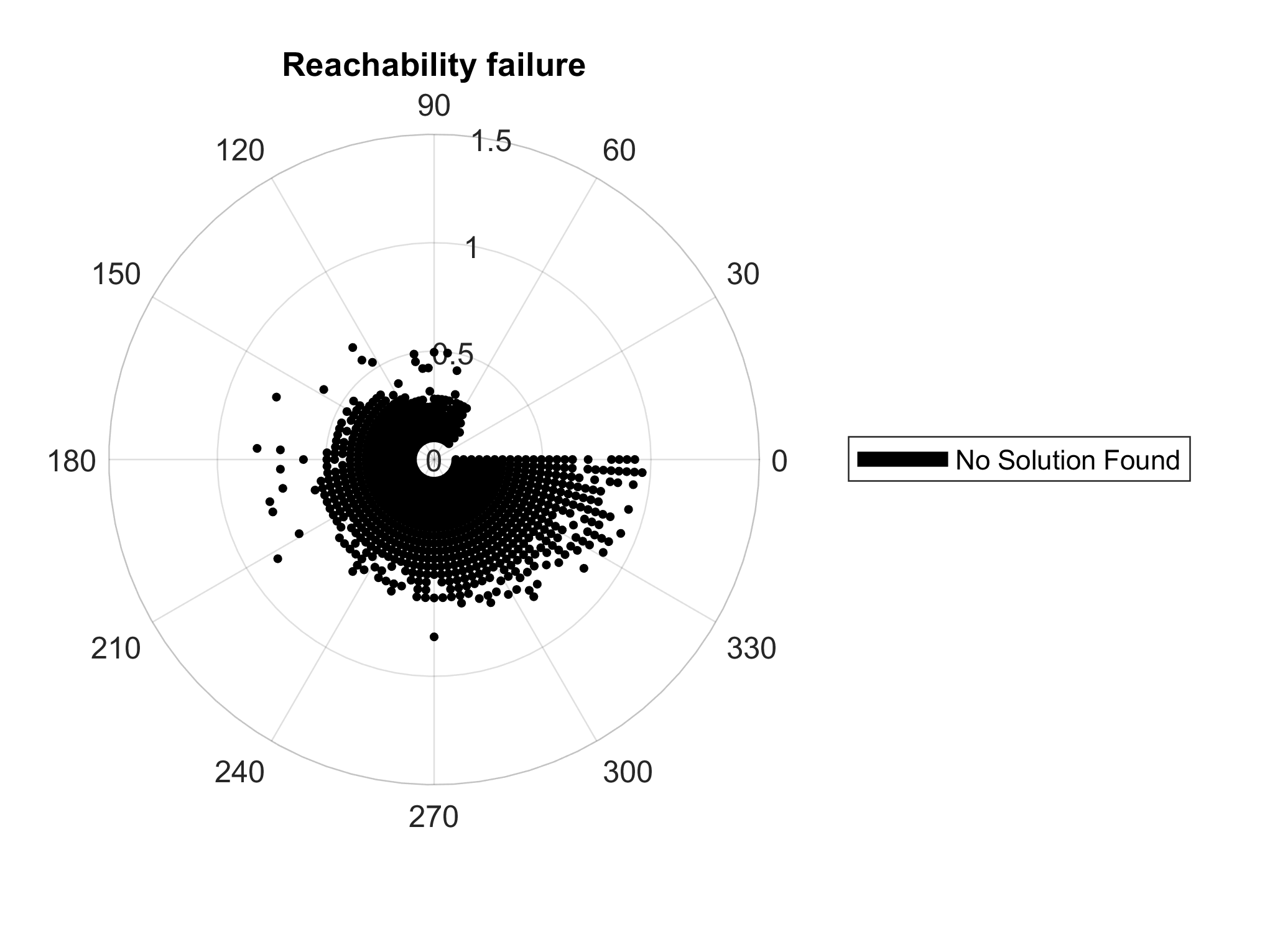}
    \caption{A modified version of Figure \ref{fig:ball_controllability}. Only points where a solution was not found are shown.}
    \label{fig:ball_failure}
\end{figure}

\section{PLANAR BOUNCING BALL}\label{sec:bouncing_ball}
As a test bed for the control problem, Problem \ref{prob:control}, we examine the controllability of a disk moving in the plane. This is a Lagrangian system with symmetry via Example \ref{ex:disk_continuous}. The Lagrangian will specifically be
\begin{equation*}
    L = \frac{1}{2}m\left(\dot{x}^2+\dot{y}^2\right) + \frac{1}{2}I\dot\theta^2 - mgy,
\end{equation*}
where $m$ is its mass, $I$ its moment of inertia, $g$ is the acceleration due to gravity, $(x,y)$ is the cartesean coordiantes of its center and $\theta$ is its orientation angle.
Recall that this Lagrangian is invariant under rotations and thus the mechanical connection will be preserved under the continuous dynamics, $\mathcal{A}=\dot\theta=\mathrm{const}$.

The control on the table manifests as both its orientation and height, $(u,v)\in\mathcal{U} = S^1\times\mathbb{R}$. See Figure \ref{fig:tilted_table} for the full schematic.

The control problem is to determine a sequences of angles and heights of the table to produce a desired change in the orientation angle of the disk. The continuous dynamics generated by the Lagrangian are
\begin{equation*}
    \begin{split}
        m\ddot{x} &= 0, \\
        m\ddot{y} &= -mg, \\
        I\ddot\theta &= 0.
    \end{split}
\end{equation*}
To determine the impact law, the location of the impact and the constraint need to be specified. If $R>0$ is the radius of the disk, the disk strikes the surface of the table when
\begin{equation*}
    h(x,y,\theta;u,v) = x\sin u(t) - y\cos u(t) + R + v(t) = 0.
\end{equation*}
Differentiating yields
\begin{equation*}
    \begin{split}
        dh &= \sin u(t) dx - \cos u(t) dy \\
        &\qquad + \left[ \left(x\cos u(t) + y\sin u(t)\right)\dot u(t) + \dot{v} \right] dt
    \end{split}
\end{equation*}
The impact constraint of rolling without slipping is 
\begin{equation*}
    \tilde{\eta} = Rd\theta + (\cos u) dx + (\sin u) dy.
\end{equation*}
The impact law is
\begin{equation*}
    \begin{split}
        m\dot{x}^+ - m\dot{x}^- &= \varepsilon\sin u + \lambda\cos u, \\
        m\dot{y}^+ - m\dot{y}^- &= -\varepsilon\cos u + \lambda\sin u, \\
        I\dot\theta^+ - I\dot\theta^- &= \lambda R, \\
        -E_L^+ + E_L^- &= \varepsilon \left[ \left(x\cos u(t) + y\sin u(t)\right)\dot u(t) + \dot{v} \right] \\
        0 &= R\dot\theta^+ + \dot{x}^+\cos u + \dot{y}^+\sin u.
    \end{split}
\end{equation*}
Under the simplification that $u=0$, solving for the two multipliers, $\varepsilon$ and $\lambda$, yields 
\begin{equation*}
    \begin{split}
        p_x^+ &= \frac{Rm}{I+mR^2} \left( Rp_x-p_\theta\right), \\
        p_y^+ &= \sqrt{ p_x^2 + [p_y-m(\dot u x+\dot{v})]^2 + \frac{m}{I}p_\theta^2 - \frac{m}{I+mR^2}(Rp_x-p_\theta)^2 } \\
        & \qquad + m(\dot u x+\dot{v}),\\
        p_\theta^+ &= \frac{I}{I+mR^2}\left( p_\theta - Rp_x\right),
    \end{split}
\end{equation*}
where $p_x = m\dot{x}$, $p_y = m\dot{y}$, $p_\theta = I\dot\theta$ and the superscripts on the right side are omitted. Notice that the angular velocity is \underline{not} conserved.

For the purposes of simplicity, we will assume that the tabletop is stationary at the moment of impact, $\dot{u} = \dot{v} = 0$. Under this simplification, the general reset map for $u\ne 0$ can be recovered via a rotational change of variables.

\subsection{Numerical Results}\label{sec:numerics}





The control problem for the bouncing ball is to find table orientations and heights such that for a specified final time $T>0$, 
\begin{equation*}
    \begin{array}{ll}
        x(0) = x_0, & x(T) = x_0, \\
        y(0) = y_0, & y(T) = y_0, \\
        \theta(0) = \theta_0, & \theta(T) = \theta_f, \\
    \end{array}
\end{equation*}
and the ball starts and ends at rest. Due to the rotational symmetry of the problem, the problem does not depend on the particular initial value of $\theta$, thus we are only interested in obtaining a desired change of angle. For concreteness, we set the initial/final conditions to $x_0=0$, $y_0=1$, and $\theta_0=0$. Now the controllability problem only depends on the final time, $T$, and the final angle, $\theta_f$.

\begin{figure}
    \centering
    \includegraphics[width=0.45\textwidth]{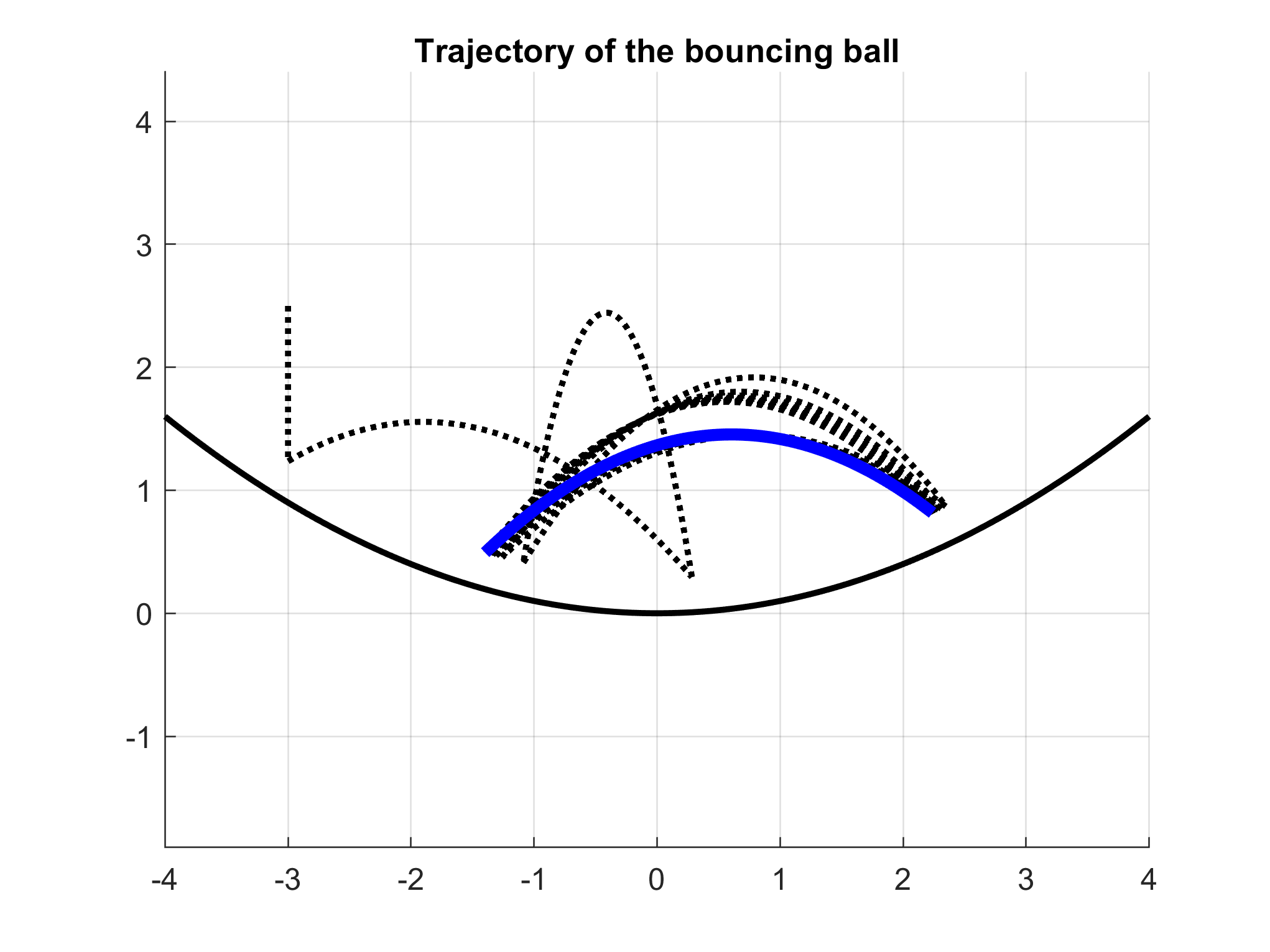}
    \caption{The ball starts at rest and its trajectory is eventually asymptotic to a periodic cycle shown in blue.}
    \label{fig:parabola_limit}
\end{figure}

\begin{figure}
    \centering
    \includegraphics[width=0.45\textwidth]{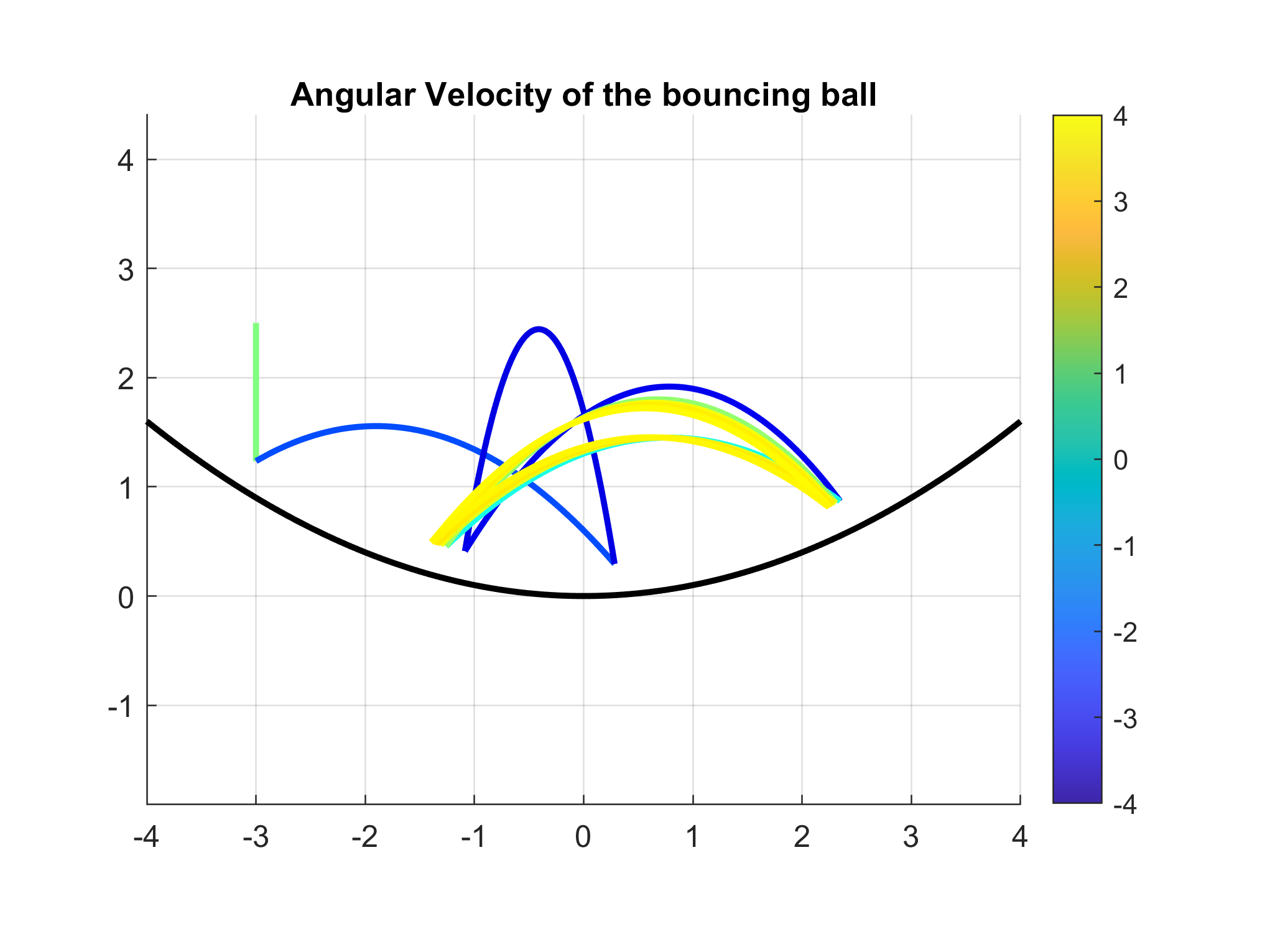}
    \caption{The same trajectory shown in Figure \ref{fig:parabola_limit}. However, here the angular velocity of the ball is recorded; its color depicts its angular velocity (rad/sec). Interestingly, the ball begins rotating clockwise while its steady-state behavior is counter-clockwise rotation.}
    \label{fig:parabola_velocity}
\end{figure}

\begin{table}[h]
    \vspace{2ex}
    \centering
    \renewcommand{\arraystretch}{2}
    \begin{tabular}{| c |c |c || c | c | }
    \hline
     \multicolumn{5}{|c|}{Parameters} \\
     \hline
     \mbox{} & Range & Resolution & Iterations & Error \\
     \hline
     $\theta_f$ & $[0, 2\pi]$ & 100 & \multirow{2}{*}{10} & \multirow{2}{*}{$\leq 10^{-2}$}\\
     \cline{1-3}
     $T$ & $[0.1, 1.5]$ & 40 &  & \\
    \hline
    \end{tabular}
    \bigskip

    \caption{
    The specific parameters used in the simulation for Figures \ref{fig:ball_controllability} and \ref{fig:ball_failure}; range refers to the range of parameters graphed; resolution refers to the amount of parameters chosen within the range.}
    \label{tab:parameters}
\vspace{-4ex}
\end{table}

This control problem was simulated via MATLAB, with the results shown in Figures \ref{fig:ball_controllability} and \ref{fig:ball_failure}. The simulation is set to allow a maximum of $5$ bounces, for which there are $10$ parameters in the initial guess: $5$ parameters for the height of the table at each bounce chosen randomly in a range of $[0,1]$, and $5$ parameters for the angle of the table, chosen randomly between $[0, \pi]$. For every data point, the simulation undergoes some number of iterations to solve the problem by randomly choosing initial guesses within the specified ranges. If there are multiple solutions with an error less than $\varepsilon >0$, the solution requiring the least number of bounces is chosen. 
The error is determined by the following equation: 
\begin{equation*}
    E = \langle q,Aq\rangle, \quad q = 
    \begin{bmatrix}
        (x_f-x_0)^2 \\ (y_f-y_0)^2 \\ (\theta_f-\theta_0)^2 \\ (T_f - T)^2
    \end{bmatrix}, \quad 
    A = \begin{bmatrix}
        1 & 0 & 0 & 0 \\
        0 & 1 & 0 & 0 \\
        0 & 0 & 50 & 0 \\
        0 & 0 & 0 & 5
    \end{bmatrix}.
\end{equation*}
Each point on the graph is defined by the final desired angle $\theta$, the radius which represents the time frame, and the color, which represents the minimum number of bounces found. The specific details are outlined in Table \ref{tab:parameters}. Note that we have chosen not to add any weight to the final angular momentum, as its magnitude was relatively insignificant in these scenarios. 

Both Figures \ref{fig:ball_controllability} and \ref{fig:ball_failure} show three features that warrant discussion: a zero bounce wedge, a lack of symmetry, and a spiral pattern. These features can be attributed to specific design choices when creating the simulation.

The size of the zero bounce wedge can be correlated with the size of the error and the weights in the error function. For very small time increments and angles, simply doing nothing gives a small enough error for a "viable" solution to be found. 


The abrupt edge and the lack of symmetry is a product of the way the angle is defined in the simulation; $\theta = 0$ is not identified with $\theta = 2\pi$. If this were changed, the zero bounce wedge would likely increase to accommodate the region in which the angle is small enough to result in no error, and the ``larger'' angles would be easier to reach.

In Figure \ref{fig:ball_controllability}, we notice the appearance of a two-bounce spiral cutting in the middle of the one-bounce spiral where there is a large angle, yet relatively short amount of time. Within that time frame, it is not possible to create a large enough angular momentum from a single bounce for the ball to reach the desired angle, so two bounces are needed to reach the desired outcome. For longer time frames, the ball has more time to rotate in the air, thus explaining the appearance of one-bounce solutions.

\color{black}



\section{CONCLUSIONS}\label{sec:conclusions}
We investigated the controllabilty problem for natural Lagrangian systems with symmetry subject to symmetry breaking controlled impacts. A case study of orientation control of a planar bouncing ball was presented. An obvious and immediate research extension is to determine a systematic procedure to guarantee controllability of this general class of control problems.

One potential approach to the controllability problem is to restrict the dynamics to the guard. In this case, controllability is related to the degree of nonholonomy of the impact distribution. A requisite for this technique to be valid is for the hybrid dynamics to converge to the restricted nonholonomic dynamics in a ``small bounce limit'' which is similar to strategy utilized in \cite{ames2006} to resolve Zeno states.

Independent of controllability, the dynamics generated by \eqref{eq:new_impact} are strange and do not share many properties with the slippery impact law \eqref{eq:moving_WE_condition}. First off, the former is not volume-preserving, unlike the latter. As shown in \cite{2101.11128}, volume-preserving impact systems possess almost no Zeno trajectories. It remains unclear whether or not the hybrid dynamics laid out here will be Zeno. Another qualitative difference is that exponentially stable periodic orbits are impossible under systems with the slippery impact law. This is no longer the case for impacts of the form \eqref{eq:new_impact}; the bouncing ball on a parabolic appears to have an exponentially stable limit cycle, Figure \ref{fig:parabola_limit}. 


A final point for future work would be to understand optimal control of intermittent contacts problems \eqref{eq:reset_control_hybrid}, with the controlled impact law \eqref{eq:controlled_impact}. This would be a natural extension of the work done in \cite{2111.11645}.

\addtolength{\textheight}{-12cm}   





\section*{ACKNOWLEDGMENT}

The authors would like to thank Dr. Jessy Grizzle for insightful conversations on this topic.

\newpage


\bibliographystyle{ieeetr}
\bibliography{references}
\end{document}